\documentclass[dvips]{amsart}
\usepackage{math,lamsarrow,pb-diagram,pb-lams,epsfig}
\newcommand\lpresn{{$L$-pre\-sen\-ta\-tion}}
\newcommand\lpresd{{$L$-pre\-sen\-ted}}
\newcommand\HNN{\textsf{HNN}}
\newcommand\coker{{\operatorname{coker}}}
\newcommand\SPRES[2]{\langle#1|\,#2\rangle}
\newcommand\SLPRES[4]{\SPRES{#1}{#2|\,#3|\,#4}}
\newcommand\SELPRES[3]{\SPRES{#1\|\,#2}{#3}}
\newcommand\PRES[2]{\big\langle#1\big|\,#2\big\rangle}
\newcommand\LPRES[4]{\PRES{#1}{#2\big|\,#3\big|\,#4}}
\newcommand\ELPRES[3]{\PRES{#1\big\|\,#2}{#3}}

\begin{document}
\title{Endomorphic Presentations of Branch Groups}
\author{Laurent Bartholdi}
\date\today
\email{laurent@math.berkeley.edu}
\address{\parbox{.4\linewidth}{Department of Mathematics\\
    University of California\\
    94720 Berkeley\\
    U. S. A.}}
\thanks{The author acknowledges support from the ``Swiss Mathematical
  Society'', the Hebrew University of Jerusalem, and the University of
  California at Berkeley}
\begin{abstract}
  We introduce ``endomorphic presentations'', or \lpresn s: group
  presentations whose relations are iterated under a set of
  substitutions on the generating set, and show that a broad class of
  groups acting on rooted trees admit explicitly constructible finite
  \lpresn s, generalising results by Igor Lysionok and Said Sidki.
\end{abstract}
\maketitle

\section{Introduction}
In the early 80's, Rostislav Grigorchuk defined a group, $G$, endowed
with many interesting properties: it is a finitely generated,
infinite, torsion group; it has intermediate growth; it has a solvable
word problem; it has finite width; etc. There are connections of $G$
to innumerably many branches of mathematics: random walks on graphs,
Hecke operators, classification of finite-rank Lie algebras,
cryptography, etc.

Already in his early papers~\cite{grigorchuk:gdegree}, Rostislav
Grigorchuk showed that $G$ is not finitely presentable. However, Igor
Lysionok obtained in~\cite{lysionok:pres} a recursively defined,
infinite set of relators for $G$, obtained by iterating a simple
letter substitution on a finite set of relators (see
Theorem~\ref{thm:G}):
\begin{quote}
  The Grigorchuk group $G$ admits the following presentation:
  \[G = \PRES{a,c,d}{\sigma^i(a^2),\sigma^i(ad)^4,
    \sigma^i(adacac)^4\;(i\ge0)},\]
  where $\sigma:\{a,c,d\}^*\to\{a,c,d\}^*$ is defined by
  $\sigma(a)=aca,\sigma(c)=cd,\sigma(d)=c$.
\end{quote}

Rostislav Grigorchuk then used this result to construct a finitely
presented amenable but not elementary-amenable
group~\cite{grigorchuk:amenEG}, thus answering negatively to an old
question by Mahlon Day~\cite{day:amen}: ``can every amenable finitely
presented group be obtained from finite groups using exact sequences
and unions?''  He also proved the independence of the relators, by
explicitly computing the Schur multiplier $H_2(G,\Z)$
in~\cite{grigorchuk:bath}.

A (friendly) competitor of the Grigorchuk group is the Gupta-Sidki
group $\doverline\Gamma$, which is also a finitely generated infinite
torsion group of subexponential growth. Both groups share other
properties, as well --- see for
instance~\cite{bartholdi-g:cras-spectrum,bartholdi-g:parabolic} where they
are studied simultaneously. Said Sidki described
in~\cite{sidki:pres} a general method yielding recursive
presentations of such groups, and for $\doverline\Gamma$ derived an
explicit, if somewhat lengthy, presentation.

Narain Gupta also followed a completely different path in obtaining
recursively presented torsion groups --- namely, he started by
defining a presentation, that is recursive but presents no explicit
regularity like the presentations considered here, and then proves
that the associated group is infinite, torsion and finitely
generated~\cite{gupta:recursive}.

In this paper, we define a general class of presentations, which we
call \emdef{endomorphic} or \emdef{$L$-presentations}. As a first
approximation, they are given by a generating set, some initial
relations, and word substitution rules that produce more relations.

We start by deriving some of their properties, and give explicit
presentations for $\doverline\Gamma$ and other contracting branch
groups (see Definitions~\ref{defn:branch} and~\ref{defn:contract}; the
main property of a \emph{branch group} is that it contains a subgroup
$K$ containing a copy $K_1\cong K^d$ for some $d\ge2$, all inclusions
having finite index. It is \emph{contracting} if there is a metric on
$K$ contracted, up to an additive constant, by the projections
$K_1\to 1^i\times K\times1^{d-i-1}$). Our main result on groups acting on rooted
trees is the following (see Theorems~\ref{thm:branch} and~\ref{thm:notfp}):
\begin{thm}
  Let $G$ be a finitely generated, contracting, semi-fractal, regular
  branch group. Then $G$ is finitely \lpresd. However, $G$ is not
  finitely presented.
  
  The Schur multiplier of $G$ has the form $A\oplus B^\infty$ for
  finite-rank groups $A,B$.
\end{thm}

\begin{defn}
  An \emdef\lpresn\ is an expression of the form
  \[L = \LPRES SQ\Phi R,\]
  where $S$ is an alphabet (i.e.\ a set of symbols), $Q,R\subset F_S$
  are sets of reduced words (where $F_S$ is the free group on $S$),
  and $\Phi$ is a set of free group homomorphisms $\phi:F_S\to F_S$.
  
  $L$ is \emdef{finite} if $S,Q,\Phi,R$ are finite. It is
  \emdef{ascending} if $Q$ is empty.

  $L$ gives rise to a group $G_L$ defined as
  \[G_L = F_S \Big/\Big\langle Q\cup
  \bigcup_{\phi\in\Phi^*}\phi(R)\Big\rangle^\#,\] where
  $\langle\cdot\rangle^\#$ denotes normal closure and $\Phi^*$ is the
  monoid generated by $\Phi$, i.e.\ the closure of $\{1\}\cup\Phi$
  under composition.

  As is customary, we shall identify the presentation $L$ and the
  group $G_L$ it defines, and write $G$ for both.
\end{defn}
The name ``\lpresn'' comes both as a homage to Igor Lysionok who
discovered such a presentation for the Grigorchuk group $G$ (see
Theorem~\ref{thm:G}) and as a reference to ``$L$-systems'' as defined
by Aristid Lindenmayer~\cite{lindenmayer:l} in the early 70's
(see~\cite{rozenberg-s:l}), used to model biological growth phenomena.

\subsection{Symmetric groups} The purpose of \lpresn s is to encode in
homomorphisms $\phi\in\Phi$ some regularity of the presentation.
Consider for instance the presentations of finite symmetric groups.
It is well-known that the following is a presentation of $\sym n$, the
symmetric group on $n$ objects
(see~\cite{burnside:presentation,moore:presentation} for its first
occurrences in literature and~\cite{sergiescu:presentations} for other
presentations):
\begin{multline*}
  \sym n=\PRES{\sigma_1,\dots,\sigma_{n-1}}{%
    \sigma_i^2\,\text{ whenever }1\le i\le n-1,\\
    (\sigma_i\sigma_{i+1})^3\text{ whenever }1\le i\le n-2,\;
    (\sigma_i\sigma_j)^2\text{ whenever }1\le i<j-1\le n-2}.
\end{multline*}

A shorter ascending \lpresn\ with the same generators can be obtained
if one lets the symmetric group act on itself by conjugation; there
remain only $3$ orbits of relators under this action. To the point,
consider the set $P=\{(1,\dots,n),(1,2),(3,\dots,n)\}$ generating
$\sym n$. For each $p\in P$, let it act as $\phi_p$ on the free group
$F_{\sigma_1,\dots,\sigma_{n-1}}$ in such a way that this action is a
lift of the action of $\sym n$ by conjugation on itself, and such that
if $\sigma_i^p=\sigma_j$ then $\phi_p(\sigma_i)=\sigma_j$ --- a simple
way of selecting such a $\phi_p$ is to pick for each $\sigma_i$ a word
$W$ over $\{\sigma_1,\dots,\sigma_{n-1}\}$ of minimal length
representing $\sigma_i^p$, and setting $\phi_p(\sigma_i)=W$,
extended by concatenation. We then obtain
\[\sym n=\ELPRES{\sigma_1,\dots,\sigma_{n-1}}{\{\phi_p\}_{p\in P}}{%
  \sigma_1^2,(\sigma_1\sigma_2)^3,(\sigma_1\sigma_3)^2}.\]
Indeed all relations $\sigma_i^2$ and $(\sigma_i\sigma_{i+1})^3$ are
obtained as $\phi_{(1,\dots,n)}^{i-1}(\sigma_1^2)$ and
$\phi_{(1,\dots,n)}^{i-1}((\sigma_1\sigma_2)^3)$, and all relations
$(\sigma_i\sigma_j)^2$ are obtained as
$\phi_{(1,\dots,n)}^{i-1}\phi_{(3,\dots,n)}^{j-i-2}((\sigma_1\sigma_3)^2)$.
Conversely, all $\phi(r)$ are relations for $\phi\in \{\phi_p\}^*$ and
$r$ a relation, since the $\phi$ are endomorphisms.

Using the same reasoning, one can obtain an ascending \lpresn\ of $\sym
n$ with only two relators, if one allows more generators:
\[\sym n=\ELPRES{\sigma_{1,2},\sigma_{1,3},\dots,\sigma_{n-1,n}}{%
  \phi_{(1,2)},\phi_{(1,\dots,n)}}{%
  \sigma_{1,2}^2,\sigma_{1,2}\sigma_{2,3}\sigma_{1,3}\sigma_{2,3}},\]
where $\sigma_{i,j}$ should be interpreted as the transposition
$(i,j)$, and $\phi_p(\sigma_{i,j})=\sigma_{i^p,j^p}$.

(This regularity in the presentation is reflected by the fact that
$H_2(\sym n,\Z)=\Z/2$ is very small --- see Subsection~\ref{subs:schur}.)

\begin{problem}
  Does there exist a bound $A$ such that all symmetric groups can be
  defined by an ascending $L$-presentation $\SELPRES S\Phi R$ of total
  length $|S|+|\Phi|+|R|<A$?
\end{problem}

\subsection{Other examples} Another example is given by presentations
of the free abelian group $\Z^n$:
\[\Z^n = \PRES{x_1,\dots,x_n}{[x_i,x_j]\,\forall i,j\in\{1,\dots,n\}}.\]
It can be expressed with fewer relators as
\[\Z^n = \ELPRES{x_1,\dots,x_n}{\phi_1,\phi_2}{[x_1,x_2]},\]
with $\phi_1: x_1\mapsto x_2,x_2\mapsto x_3,\dots,x_n\mapsto x_1$ and
$\phi_2: x_1\mapsto x_1,\;x_2\mapsto x_3,x_3\mapsto
x_4,\dots,x_n\mapsto x_2$.  Of course, the main interest of \lpresn s
is to encode groups that don't even have a finite presentation:
consider for instance the group $\sym\infty\Z$ of permutations of $\Z$
that act like a translation outside a finite interval. It is generated
by $\sigma=(1,2)$ and $\tau:n\mapsto n+1$:
\begin{align*}
  \sym\infty\Z &= \PRES{\tau,\sigma}{\sigma^2,
  [\sigma,\tau^n]^2\,\forall n\ge2,[\sigma,\tau]^3}\\
  &= \LPRES{\tau,\sigma,\overline\sigma}{\sigma\overline\sigma}\phi
  {\sigma^2,[\sigma,\tau]^3,[\sigma,{\overline\sigma}^{\tau^2}]}\\
  &= \LPRES{\tau,\sigma,\overline\tau}{\tau^{-1}\overline\tau}\psi
  {\sigma^2,[\sigma,\tau]^3,[\sigma,\sigma^{\tau\overline\tau}]}
\end{align*}
with $\phi(\overline\sigma)={\overline\sigma}^\tau$ and
$\psi(\overline\tau)=\tau\overline\tau$, both preserving the other
generators $\sigma$ and $\tau$. The extra generators $\overline\sigma$
and $\overline\tau$ are just convenient copies of the generators.

\subsection{Outline} The paper is organized as follows:
Section~\ref{sec:gt} contains group-theoretical results on \lpresn s.
Section~\ref{sec:branch} describes the main result of this paper,
namely that all finitely generated regular branch groups have a finite
\lpresn. Section~\ref{sec:examples} describe the \lpresn s of the $5$
``testbed'' branch groups introduced in~\cite{bartholdi-g:cras-spectrum}.

I wish to thank, in reverse alphabetical order, Zoran \v Suni\'k,
Rostislav Grigorchuk, Denis Osin, Gulnara Arjantseva for their
entertaining discussions on this question, and their careful reading
of the text. Gulnara generously offered contributions to
Subsection~\ref{subs:embeddings}, Denis to
Subsection~\ref{subs:identities} and Zoran to
Subsection~\ref{subs:perfect}.

Many results were obtained using the software system
\textsc{Gap}~\cite{gap:manual}, whom I thank for his patient and
silent permutation-grinding.

\subsection{Notations} For me, $g^h$ denotes $h^{-1}gh$, and the
expression $g^{\sum n_ih_i}$ means $\prod h_i^{-1}g^{n_i}h_i$. The
commutator $[g,h]$ is $g^{-1}h^{-1}gh$, and $X^*$ is the monoid
generated by $X$. The normal closure of $X$ in $G$ is written $\langle
X\rangle^\#$, the normal subgroup of $G$ \emph{normally generated} by
$X$.

\section{Group-Theoretical Properties}\label{sec:gt}
In this section, we are interested in the following questions: which
group-theoretical constructions preserve the property of having a
finite \lpresn? which groups admit a finite \lpresn?

We shall say a group is \emdef{finitely \lpresd} if it admits a
finite \lpresn.
\begin{rem}
  There are finitely \lpresd\ groups that, for some imposed
  generating set, do not admit a finite \lpresn. This is in
  contrast with finitely presented groups, for which admitting a
  finite presentation is independent of the choice of generators ---
  that property is even invariant under quasi-isometries.

  For instance, consider the ``lamplighter group'' of
  Theorem~\ref{thm:ll}, with its finite \lpresn. This group
  $G$ does not have a finite \lpresn\ with generators
  $\{a,t\}$, as can be seen by a careful study of endomorphisms of
  $\F_2$.
\end{rem}

\begin{prop}
  Let $G$ admit a finite ascending \lpresn, and let $S'$ be a finite
  generating set of $G$. Then $G$ admits a finite ascending \lpresn\
  with generators $S'$.
\end{prop}
\begin{proof}
  The standard proof that being finitely presented involves Tietze
  transformations, and extends to \lpresn s. One changes $S$ into $S'$
  by a finite number of ``Tietze moves'', which either replace a
  generator by a product or quotient of generators, or add or delete a
  generator $s$ along with the relator $s$.

  For an \lpresn\ $\SELPRES S\Phi R$, the operations are as follows:
  if the move was to replace the generator $s$ by $s'=st$, one
  replaces all instances of $s$ by $s't^{-1}$ in $R$ and the images of
  $\phi\in\Phi$, modifying them by $\phi(s')=\phi(st)$.

  If the move was addition of a generator $s$ to $S$ and $R$, one
  extends all $\phi\in\Phi$ by $\phi(s)=1$. If the move was deletion
  of $s$ from $S$ and $R$, one deletes all instances of $s$ in the
  images of all $\phi\in\Phi$, and adds $\phi(s)$ to $R$.
\end{proof}

\subsection{Embeddings}\label{subs:embeddings}
We start by some motivation for the study of \lpresn s. Recall Graham
Higman's Embedding Theorem
\begin{thm}[\cite{higman:embed} or Section IV.7 of~\cite{lyndon-s:cgt}]
  A countable group $G$ can be embedded in a finitely presented group
  $\widehat G$ if and only if it is recursively presented.
\end{thm}
The first proof by Higman was unconstructive; since then, explicit
constructions of $\widehat G$ were
given~\cite{aanderaa:higman,aanderaa-c:higman,olshansky-s:higman}.
They require, however, a good mastery of Turing- or $S$-machine
programming. I am not aware of an explicit finitely presented group
containing $\mathbb Q$. In contrast, a finitely \lpresd\ group
containing $\mathbb Q$ follows: first,
\[\Q = \PRES{x_1,x_2,\dots}{x_nx_{n+1}^{-n-1}\,\forall n\ge1},\]
where $x_n$ should be interpreted as $1/n!\in\Q$. We embed $\Q$ in the
finitely generated group
\[G = \PRES{x,y}{x^ny^n(x^{n+1}y^{n+1})^{-n-1}\,\forall n\ge1}
= \PRES{x,y}{yx(x^{n+1}y^{n+1})^n\,\forall n\ge1}\]
through $x_n\mapsto x^ny^n$; now $G$ embeds in the finitely \lpresd\
group $H$ given by
\[H = \LPRES{x,y,a,b,c,d,e,c',d',e'}{c^{-1}c',d^{-1}d',e^{-1}e',[\langle
  d,e\rangle,\langle x,y\rangle]}{\phi_1,\phi_2,\phi_3,\phi_4}{%
  yxbac,(e')^{d'}c'}\]
where
\[\phi_1:\begin{cases}a\mapsto bac\end{cases}\quad
\phi_2:\begin{cases}a\mapsto1\\ b\mapsto d^{-1}ybxd\end{cases}\quad
\phi_3:\begin{cases}a\mapsto1\\ b\mapsto yex\end{cases}\quad
\phi_4:\begin{cases}c'\mapsto c'c\\ d'\mapsto d'd\\ e'\mapsto e'e\\
\end{cases},\]
it being understood that unspecified generators map to themselves.
Indeed $n-1$ applications of $\phi_1$ to $yxbac$ yield $yxb^nac^n$;
then $n$ applications of $\phi_2$ yield $yx((d^{-1}y)^nb(xd)^n)^nc^n$;
then $\phi_3$ and the commutation relations yield
$yx(y^{n+1}x^{n+1})^nd^{-n}e^nd^nc^n$. On the other hand,
$\phi_4^{n-1}((d')^{-1}e'd'c')$ yields $d^{-n}e^nd^nc^n$, so
$yx(y^{n+1}x^{n+1})^n=1$ in $H$, whence
$x^ny^n=(x^{n+1}y^{n+1})^{n+1}$ in $H$. Any other sequence of
operations $\phi_i$ would give a long relation containing
non-$\{x,y\}$ symbols, so $G$ embeds in $H$.

Some finitely \lpresd\ groups embed nicely in finitely presented
groups; recall that the \HNN\ extension
$\Omega(G,H\overset\phi\longrightarrow K)$ is \emdef{ascending} if
$H=G$.
\begin{thm}\label{thm:hnn}
  Let $G$ be finitely \lpresd\ by an ascending \lpresn.
  Then a finitely presented group $\widehat G$ containing $G$ can be
  effectively constructed.  Moreover, $\widehat G$ is an ascending
  \HNN\ extension of $G$ by a finite number of stable letters.

  In case $G$ is amenable, $\widehat G$ is a finitely presented
  amenable group containing $G$.
\end{thm}
\begin{proof}
  Let $\SELPRES S\Phi R$ be a finite ascending \lpresn\ of
  $G$. Consider the group
  \[\widehat G = \PRES{S\cup\Phi}{R
    \cup\{s^\phi=\phi(s)\}_{s\in S,\phi\in\Phi}}.\]
  It is finitely presented, and the map $G\to\widehat G$ defined by
  sending $s\in S$ to $s$ is a well-defined injective homomorphism,
  since the $\phi:S^*\to S^*$ induce injective homomorphisms of $G$.
\end{proof}

As a consequence, we may construct finitely generated subgroups of
hyperbolic groups that are not hyperbolic. Recall that a set $W$
satisfies the \emdef{small cancellation condition} $C'(\epsilon)$ if
for any $u,v\in\overline W$ the common prefix of $u$ and $v$ has
length strictly less than $\epsilon\min\{|u|,|v|\}$, where $\overline
W$ is the closure of $W$ under taking inverses and cyclic conjugates.

\begin{cor}\label{cor:nonqc}
  Let $\SELPRES S\phi R$ be a finite ascending \lpresn\ of
  $G$ with $R\neq\emptyset$, and let $\widehat G$ be the finitely
  presented group constructed above. Assume that $R\cup\phi(S)$
  satisfies the small cancellation condition $C'(1/6)$.

  Then $\widehat G$ is a hyperbolic group containing a non-hyperbolic
  finitely generated, infinitely presented subgroup $G$. In particular
  $G$ is not quasi-convex in $\widehat G$.
\end{cor}
\begin{proof}
  It follows from the hypotheses that $\bigcup_{n\ge0}\phi^n(R)$ is
  also $C'(1/6)$. Therefore $G$ is not finitely presented, so cannot
  be hyperbolic. On the contrary, $\widehat G$ is finitely presented
  and its presentation is $C'(1/6)$, so it is hyperbolic. Finally a
  quasi-convex subgroup of a hyperbolic group would be hyperbolic, so
  $G$ cannot be quasi-convex.
\end{proof}

As a simple application of this corollary, consider the group
\[G=\ELPRES{x,y}{\phi}{(xy)^7},\]
with $\phi(x)=x^7$ and $\phi(y)=y^7$, embedding in
\[\widehat G=\PRES{x,y,\phi}{x^{\phi-7},y^{\phi-7},(xy)^7}.\]

\begin{prop}
  If $G$ is finitely presented, then it is finitely \lpresd.
  There are non-finitely \lpresd\ groups, and there are finitely
  \lpresd, but not finitely presented, groups.
\end{prop}
\begin{proof}
  The first claim is obvious: finite \lpresn s with
  $R=\Phi=\emptyset$ are precisely finite presentations.

  There are only countably many finite \lpresn s, but
  uncountably many finitely-generated groups, so ``most'' groups are
  not finitely \lpresd.
  
  Finally, Theorem~\ref{thm:ll} shows that the ``lamplighter group''
  described there is finitely \lpresd, but not finitely
  presented.
\end{proof}
Note, however, that it is not trivial to explicitly point at a
non-finitely $L$-presentable group. A group having a
non-recursively-enumerable presentation satisfying some small
cancellation condition would be an example. The free group in a
variety defined by infinitely many identities (they exist
by~\cite{olshansky:fb}) is another one. More examples appear in the
course of Lemma~\ref{lem:uncountext}.

\begin{prop}
  If $G,H$ are finitely \lpresd\ groups, then $G*H$ is finitely
  \lpresd.
  
  If $G$ is finitely \lpresd\ and $H,K$ are isomorphic finitely
  generated subgroups of $G$, then the \HNN\ extension
  $\Omega(G,H\overset\psi\longrightarrow K)$ is finitely \lpresd.
\end{prop}
\begin{proof}
  Let $\SLPRES SQ\Phi R$ be a finite \lpresn\ of $G$, and let $\SLPRES
  TP\Psi U$ be a finite \lpresn\ of $H$. A finite \lpresn\ of $G*H$ is
  \[\LPRES{S\cup T}{Q\cup P}{\Phi\cup\Psi}{R\cup U},\] 
  where it is understood that each $\phi\in\Phi$ is extended to a
  homomorphism $\phi:(S\cup T)^*\to(S\cup T)^*$ by mapping each $t\in
  T$ to itself; and similarly for each $\psi\in\Psi$.
  
  Let now $H$ be the subgroup of $G$ generated by $T\subset S^*$. A
  presentation for the \HNN\ extension of $G$ by $\psi:H\to K$ is
  \[\LPRES{S\cup\{\psi\}}{Q\cup\{\psi(t)^{-1}t^\psi\}_{t\in T}}\Phi R.\]
\end{proof}

\begin{prop}\label{prop:ext}
  If $G,H$ are finitely \lpresd\ groups, then any split extension
  of $G$ by $H$ is finitely \lpresd. If $H$ is finitely
  presented, then any extension of $G$ by $H$ is finitely
  \lpresd.
\end{prop}
\begin{proof}
  Let $\SLPRES SQ\Phi R$ be a finite \lpresn\ of $G$; let $\SLPRES
  TP\Psi U$ be a finite \lpresn\ of $H$; let $X$ be an extension of
  $G$ by $H$, given as $1\to G\to X\to H\to 1$. Let $\sigma$ be a
  section of $H$ to $X$; in case the extension splits, we suppose that
  $\sigma$ is a group homomorphism.
  
  Each relator $p\in P$ lifts through $\sigma$ to an element $g_p\in
  G$, so we may define $P'=\{pg_p^{-1}|\,p\in P\}$, a set of relators
  in $X$. Since $G$ is normal in $X$, we also have
  $s^{\sigma(t)}=g_{s,t}\in G$ for each $s\in S,t\in T$. Consider now
  the presentation
  \begin{equation}\label{eq:ext}
    \LPRES{S\cup T}{Q\cup P'\cup\{s^tg_{s,t}^{-1}\}_{s\in S,t\in T}}{%
    \Phi\cup\Psi}{R\cup U},
  \end{equation}
  where it is understood that each $\phi\in\Phi$ is extended to a
  homomorphism $\phi:(S\cup T)^*\to(S\cup T)^*$ by mapping each $t\in
  T$ to itself; and similarly for each $\psi\in\Psi$.
  
  If $X$ is a split extension, then $g_p=1$ for all $p\in P$, and
  similarly all $\phi(u)$ (with $u\in U$ and $\phi\in\Phi^*$) are
  relations in $X$. If $H$ is finitely presented, we may suppose
  $U=\emptyset$ and again all relations given in~(\ref{eq:ext}) are
  satisfied.
  
  We have shown that in the cases considered $X$ is a quotient
  of~(\ref{eq:ext}). Let now $w$ be a word in $S\cup T$ equal to $1$
  in $X$. The relations $s^t=g_{s,t}$ allow $w$ to be written as
  $s_1\dots s_nt_1\dots t_m$; then projecting to $H$ gives $t_1\dots
  t_m=1$ by applying relations in $H$. The same relations
  in~(\ref{eq:ext}) will reduce $s_1\dots s_nt_1\dots t_m$ to a word
  in $S^*$; the corresponding element of $G$ can be reduced to $1$
  using relations in $G$, and these same relations exist in $X$, so
  $w=1\in X$ and~(\ref{eq:ext}) is a presentation of $X$.
\end{proof}

Note that there are extensions of finitely $L$-presented groups that
are not finitely $L$-presented; more precisely,
\begin{lem}\label{lem:uncountext}
  There are uncountably many non-isomorphic extensions of $\Z/2$ by
  $\Z/2\wr\Z$.

  As a consequence, there are uncountably many such extensions that
  are not finitely $L$-presented.
\end{lem}
\begin{proof}
  $H^2(\Z/2\wr\Z,\Z/2)=(\Z/2)^\infty$ --- see Subsection~\ref{subs:schur}.
\end{proof}

\begin{prop}
  If $G$ is an finitely \lpresd\ group, then any finite-index
  subgroup of $G$ is finitely \lpresd.
  
  If $N\triangleleft G$ is finitely generated as a normal subgroup of
  a finitely \lpresd\ group $G$, then $G/N$ is finitely \lpresd.
\end{prop}
\begin{proof}
  Let $\SLPRES SQ\Phi R$ be a finite \lpresn\ of $G$, and let $X$ be a
  right transversal of the finite-index subgroup $H$ of $G$.  In view
  of Proposition~\ref{prop:ext}, we may suppose $H$ is normal in $G$,
  since any finite-index subgroup is a finite extension of its core,
  which is normal of finite index.
  
  We then have $G=\bigcup_{x\in X}Hx=\bigcup_{x\in X}xH$. For $g\in
  G$, let $\overline g\in X$ denote its coset representative. By the
  Reidemeister-Schreier method, $H$ is generated by the finite set
  $T=\big\{s^x\big\}_{x\in X,s\in S}$, and a presentation of $H$ is
  given by
  \[\Big\langle T\Big|\,
  \big\{\widetilde{q^x}\big\}_{q\in Q,x\in X}\cup
  \big\{\widetilde{\phi(r)^x}\big\}_{\phi\in\Phi^*,r\in R,x\in
    X}\Big\rangle,\]
  where $\widetilde w$ is a rewriting of $w$ as a word over $T$.  Now
  each $\phi\in\Phi$ induces naturally a monoid homomorphism
  $\widetilde\phi$ over $T^*$, and since
  $\widetilde{\phi(r)^x}=\widetilde\phi(\widetilde{r^x})$, a finite
  \lpresn\ for $H$ is given by
  \[\Big\langle T\Big|\,
  \big\{\widetilde{q^x}\big\}_{q\in Q,x\in X}\Big|\,
  \{\widetilde\phi\}_{\phi\in\Phi}\Big|\,
  \{\widetilde{r^x}\}_{r\in R,x\in X}\Big\rangle.\]
  
  For the second statement of the proposition, let $\SLPRES SQ\Phi R$
  be a finite \lpresn\ of $G$ and let $T$ be a finite generating set
  for $N$. Then
  \[\LPRES{S}{Q\cup N}\Phi R\]
  is a finite \lpresn\ of $G/N$.
\end{proof}

\begin{prop}
  If $G,H$ are finitely \lpresd\ groups, and either $G$ is abelian or
  $H$ is finite, then the restricted wreath product $G\wr H$ is
  finitely \lpresd.
\end{prop}
\begin{problem}
  The corresponding assertion with ``finitely \lpresd'' replaced by
  ``finitely presented'' does not hold. Under which conditions does
  the statement hold, for non-abelian $G$ and infinite $H$?
\end{problem}
\begin{proof}
  If $H$ is finite, then $G\wr H$ is finitely \lpresd\ by
  Proposition~\ref{prop:ext}. Let us assume then that $G$ is abelian.
  Let $\SLPRES SQ\Phi R$ be a finite \lpresn\ of $G$, and let $\SLPRES
  TP\Psi U$ be a finite \lpresn\ of $H$. An \lpresn\ of $G\wr H$ is
  \[\LPRES{S\cup T}{Q\cup P\cup\{[s_1,s_2^h]\}_{s_1,s_2\in
  S,h\in H}}{\Phi\cup\Psi}{R\cup U},\]
  where it is understood that each $\phi\in\Phi$ is extended to a
  homomorphism $\phi:(S\cup T)^*\to(S\cup T)^*$ by mapping each $t\in
  T$ to itself; and similarly for each $\psi\in\Psi$. This
  \lpresn\ is in general not finite, but this can be remedied
  by introducing new generators $\overline S$ in bijection with $S$
  and new homomorphisms $\Omega_T$ in bijection with $T$:
  \[G\wr H=\LPRES{S\cup T\cup\overline S}{%
    Q\cup P\cup\{s^{-1}\overline s\}_{s\in S}}{\Phi\cup\Psi\cup\Omega_T}{%
    R\cup U\cup\{[s_1,\overline s_2]\}_{s_1,s_2\in S}},\] where
  $\omega_t\in\Omega_T:(S\cup T\cup\overline S)^*\to(S\cup
  T\cup\overline S)^*$ is defined by $\omega_t(\overline s)=\overline
  s^t$ and $\omega_t(s)=s$ and $\omega_t(t')=t'$.  Indeed the new
  generators $\overline s$ do not enlarge $G$, since $\overline s=s$
  is a relation; also,
  \[[s_1,s_2^h]=[s_1,\overline{s_2}^h]=[s_1,\overline{s_2}^{t_1\dots
    t_n}]=[s_1,\omega_{t_1}\dots\omega_{t_n}(\overline{s_2})]=
  \omega_{t_1}\dots\omega_{t_n}([s_1,\overline{s_2}])=1\]
  is a relation, for all $h=t_1\dots t_n\in H$.
\end{proof}

\begin{problem}
  Let $G$ be a finitely $L$-presented group generated by $S$, let $H$
  a finitely generated subgroup, and let $X$ be a transversal of $H$
  in $G$ which is a regular subset of $S^*$. Under which extra
  conditions is $H$ finitely \lpresd?
\end{problem}

\subsection{Identities}\label{subs:identities}
We now show that groups defined by identities are all finitely
\lpresd. Recall that an identity is a word $w\in F_Y$, and that the
group $G$ \emdef{satisfies the identity $w$} if $f(w)=1$ for all
$f:F_Y\to G$. For instance, all abelian groups satisfy the identity
$[y_1,y_2]$. The \emdef{free group on $S$ with respect to $w$} is
$F_S/\langle f(w)\quad\forall f:F_Y\to F_S\rangle$. It is the largest
group satisfying $w$, in the sense that every group generated by $S$
and satisfying $w$ is a quotient of it.  These groups are also
sometimes referred to as \emdef{relatively free groups} of finitely
based varieties~\cite{neumann:varieties}. In that spirit, a group has
\emdef{presentation $\SPRES SR$ within a variety} if it is the
quotient of the free group on $S$ in that variety by the normal
closure of $R$.

\begin{prop}
  Let $G$ be finitely generated and finitely presented with respect to
  the identity $w$. Then $G$ is finitely \lpresd.
\end{prop}
\begin{proof}
  It suffices to prove the claim for a relatively free group, since
  the quotient of a finitely \lpresd\ group by a finitely normally
  generated normal subgroup remains finitely \lpresd.
  
  Let us then suppose $G$ relatively free, and generated by $X$, and
  write $w=w(y_1,\dots,y_n)\in F_Y$.  For $x\in X^{\pm1}$ and $y\in
  Y$, define the endomorphism $\phi_{xy}$ of $F_{X\sqcup Y}$ by
  $\phi{xy}(y)=xy$, and $\phi_{xy}(s)=s$ for all other $s\in X\sqcup
  Y$. Then the following is a finite \lpresn\ of $G$:
  \[\LPRES{X\sqcup Y}{Y}{\{\phi_{xy}\}_{x\in X^{\pm1},y\in Y}}{\{w\}}.\]
  Indeed write $\Phi=\{\phi_{xy}\}_{x\in X^{\pm1},y\in Y}$. Then
  \begin{align*}
    \LPRES{X\sqcup Y}Y\Phi{\{w\}} &= \PRES{X\sqcup Y}{Y\cup\Phi(w)}
    = \PRES{X\sqcup Y}{Y\cup\{w(w_1(X)y_1,\dots,w_n(X)y_n)\}}\\
    &= \PRES{X}{\{w(w_1(X),\dots,w_n(X))\}} = \PRES{X}{f(w)\quad\forall f:F_Y\to F_X},
  \end{align*}
  where the $w_i(X)$ are arbitrary words over $X$, and $f:F_Y\to F_X$
  is given by $f(y_i)=w_i(X)$.
\end{proof}

As a consequence, the free Burnside groups (defined by the identity
$a^n\in F_a$), the rank-$r$ free solvable groups, etc.\ are finitely
\lpresd.  Moreover:
\begin{cor}
  Any finitely generated group in the variety $\mathfrak A\mathfrak
  N_k$ of abelian-by-(nilpotent of degree $k$) groups is finitely
  \lpresd.
\end{cor}
\begin{proof}
  By~\cite{hall:soluble}, every group in the variety $\mathfrak
  A\mathfrak N_k$ is the quotient of the free group in that variety
  (defined by the identity $[[x_1,\dots,x_k],[y_1,\dots,y_k]]\in
  F_{x_i,y_i}$) by a finite number of relations.
\end{proof}

\subsection{Schur multipliers}\label{subs:schur}
It is well known, by Issai Schur and Heinz Hopf's
formula~\cite[Theorem~5.3]{brown:cohomology}, that the Schur
multiplier $H_2(G,\Z)$ ($=H^2(G,\C^\times)$ for finite groups) of a
group $G$ can be computed from a presentation of $G$; namely, given a
presentation $G=\SPRES ST$, we have
\[H_2(G,\Z)=\frac{\langle T\rangle^\#\cap[F_S,F_S]}{[\langle T\rangle^\#,F_S]}.\]
As a consequence, a finitely presented group necessarily has a
finite-rank Schur multiplier. (Note, however, that the converse is not
true --- see Theorem~\ref{thm:baumslag}.) We note that Hopf's formula
extends to \lpresn s.

Let us first recall a few facts on Schur multipliers;
see~\cite{karpilovsky:schur} for further details:
\begin{itemize}
\item Norman Blackburn's result~\cite{blackburn:wreath}
  \begin{equation}\label{eq:blackburn}
    H_2(H\wr G,\Z)=H_2(G,\Z)\oplus H_2(H,\Z)\oplus\frac{\{f:G\to H/H'\otimes
    H/H'\}}{\{f(x^{-1})=\tau f(x)\}},
  \end{equation}
  where $\tau:H/H'\otimes H/H'\to H/H'\otimes H/H'$ sends $h\otimes
  h'$ to $h'\otimes h$, and the $f$ above are just set maps.
\item A special case of the K\"unneth's formula,
  \begin{equation}\label{eq:product}
    H_2(G\times H,\Z)=H_2(G)\oplus H_2(H)\oplus(G/G'\otimes H/H'),
  \end{equation}
\item Shapiro's lemma: for an exact sequence $1\to N\to G\to Q\to1$,
  \begin{equation}\label{eq:shapiro}
    H_2(N,\Z)=H_2(G,\Z Q),
  \end{equation}
  with the $G$-action on $\Z Q$ induced by the quotient map $G\to Q$.
\end{itemize}

\begin{thm}
  Let $G$ admit a finite \lpresn\ $\SLPRES SQ\Phi R$.
  Then $H_2(G,\Z)=A\oplus \bigoplus_{\Phi^*}B$, where $A$ and $B$ are
  finitely-generated abelian groups.
\end{thm}

\begin{proof}
  Write $F=F_S$, and $W=\langle\Phi^*(R)\rangle^\#$.
  Consider the group $W/[W,F]$. It is abelian and generated by
  $\Phi^*(R)$. The maps $\phi\in\Phi$ are such that $\coker\phi$ is
  finitely generated (by $R$), and we may filter $R$ along
  $\Phi^*$. For each $\phi\in\Phi$, write
  \[1\longrightarrow\ker\phi\longrightarrow
  W\overset\phi\longrightarrow W\longrightarrow A_\phi\oplus
  B_\phi\longrightarrow1,\]
  where $A_\phi$ splits back into $W$ and $B_\phi$ does not. Then
  $W$ lies inside
  $\bigoplus_{\phi\in\Phi}B_\phi\oplus\bigoplus_{\phi\in\Phi^*}A_\phi$,
  so $W$ is of the required form. The Schur multiplier is obtained
  from $W$ by extending by $\langle S\rangle^\#/[F,S]$ (which has
  finite rank), and restricting to $[F,F]$, both operations preserving
  the claimed form of $H_2(G,\Z)$.
\end{proof}

It follows, for instance, that $H_2(G,\Z)$ may be neither $\Q$ nor
$\Z[\frac1n]$, for a finitely \lpresd\ group. However, the Schur
multiplier may be trivial, as in Theorem~\ref{thm:baumslag}, or
of infinite-rank, as in the examples of branch groups of
Subsection~\ref{subs:br:schur}.

\section{Branch Groups}\label{sec:branch}
The purpose of this section is to prove the following general results:
\begin{thm}\label{thm:branch}
  Let $G$ be a finitely generated, contracting, semi-fractal, regular
  branch group. Then $G$ is finitely \lpresd.
\end{thm}

\begin{thm}\label{thm:notfp}
  Let $G$ be a finitely generated, contracting, semi-fractal, regular
  branch group. Then $G$ is not finitely presented.
\end{thm}
Even though I am sure that the contracting hypothesis is not needed in
Theorem~\ref{thm:notfp}, I have been unable to prove it without that
extra condition.

We start by recalling some definitions
from~\cite{bartholdi-g:cras-spectrum,grigorchuk:jibg} concerning
branch groups. We fix an integer $d\ge2$, and the finite alphabet
$\Sigma=\Z/d\Z$, written $\{1,\dots,d\}$ for convenience.  The
\emdef{$d$-regular rooted tree} is the free monoid $\Sigma^*$.  A tree
automorphism $g\in\aut\Sigma^*$ is a bijective map
$g:\Sigma^*\to\Sigma^*$ that preserves prefixes, i.e.\ such that
$g(\sigma\tau)\in g(\sigma)\Sigma^{|\tau|}$ for all
$\sigma,\tau\in\Sigma^*$. There is an isomorphism between the subtree
$\sigma\Sigma^*$ and $\Sigma^*$, given by left-cancellation of
$\sigma$. It induces an isomorphism
$\pi_\sigma:\aut(\sigma\Sigma^*)\to\aut(\Sigma^*)$.

A \emdef{$d$-rooted group} is a finitely generated subgroup $G$ of
$\aut\Sigma^*$. The \emdef{rooted automorphism} is the automorphism
$a\in\Sigma^*$ defined by
\[a(\sigma_1\sigma_2\dots\sigma_n) =
(\sigma_1+1)\sigma_2\dots\sigma_n.\]

Fix a rooted group $G$, let $\stab_G(\sigma)$ be the stabilizer in $G$
of the vertex $\sigma\in\Sigma^*$, and set
$\stab_G(n)=\bigcap_{\sigma\in\Sigma^n}\stab_G(\sigma)$. Restriction
induces a map
\[\pi_\sigma:\stab_G(\sigma)\to\aut(\sigma\Sigma^*)\to\aut\Sigma^*.\]
The group $G$ is \emdef{fractal} if $\pi_\sigma\stab_G(\sigma)=G$ for
all $\sigma\in\Sigma^*$, and \emdef{semi-fractal} if
$\pi_\sigma\stab_G(\sigma)\le G$ for all $\sigma\in\Sigma^*$. In that
case, the map
\[\psi=(\psi_1,\dots,\psi_d)\omega:\stab_G(1)\to G^\Sigma\]
defined by $\psi_i(g)=\pi_i(g_{|i\Sigma^*})$ is an embedding. It
extends to a map still written $\psi:G\to G\wr\sym\Sigma$, by lifting
$\psi$ to $G$ using the natural map $G\to\sym\Sigma$ given by
restriction to the first level of the tree.

The \emdef{rigid stabilizer} of the vertex $\sigma$ is
\[\rist_G(\sigma)=\bigcap_{\tau\not\in\sigma\Sigma^*}\stab_G(\tau)\]
and the \emdef{rigid level stabilizer} of level $n$ is
\[\rist_G(n) = \prod_{\sigma\in\Sigma^n}\rist_G(\sigma).\]
Note $\rist_G(\sigma)<\stab_G(\sigma)$ and $\rist_G(n)<G$ for all
$\sigma\in\Sigma^*$ and $n\in\N$.

The group $G$ is \emdef{level-transitive} if it acts transitively on
$\Sigma^n$ for all $n\in\N$. In that case, $\stab_G(\sigma)$ and
$\rist_G(\sigma)$ depend, up to isomorphism, only on the length of
$\sigma$.

\begin{defn}\label{defn:branch}
  The group $G$ is a \emdef{branch} group if it is level-transitive,
  and $\rist_G(n)$ has finite index in $G$ for all $n$. It is
  \emdef{weak branch} if all $\rist_G(\sigma)$ are non-trivial (and
  hence infinite).  It is \emdef{regular branch} if
  $[G:\pi_\sigma\rist_G(\sigma)]$ is (finite and) constant for all
  long enough $\sigma\in\Sigma^*$. In that case, there is a
  finite-index subgroup $K\le G$ such that $K^\Sigma\le\psi(K)$, and
  $G$ is \emdef{regular branch over $K$}.
\end{defn}

\begin{defn}\label{defn:contract}
  Let $G$ be a branch group generated by a finite set $S$, and consider
  the induced word metric on $G$. We say $G$ is \emdef{contracting} if
  there is a constant $D$ such that for every word $w\in S^*$
  representing an element of $\stab_G(1)$, writing
  $\psi(w)=(w_1,\dots,w_d)$, we have
  \begin{equation}\label{eq:contract}
    |w_i|<|w|\text{ for all }i\in\Sigma,\text{ as soon as }|w|>D.
  \end{equation}
\end{defn}
It then follows that there is an algorithm $\mathcal A$ solving the
word problem in $G$: in this algorithm, we only assume that given a
group generator we know its action on the top level $\Sigma$ of the
tree, and that given a word representing an element of $\stab_G(1)$ we
may compute explicitly $\psi(w)$.

\begin{description}
\item[Initialization] Let $V\subset S^*$ be the set of all words of
  length at most $D$, and let $W\subset V$ be the set of words acting
  trivially on $\Sigma$. Note that $\psi$ is a well-defined map from
  $W$ to $V^d$. Assign to each $v\in V$ a flag, that is either
  ``trivial'', ``non-trivial'' or ``unknown yet''.  Initially all
  flags are ``unknown yet''.
  
  For each $v\in V$ flagged ``unknown yet'', if $v\in V\setminus W$ or
  $\psi(v)$ has a component flagged ``non-trivial'', flag $v$ as
  ``non-trivial''. If $\psi(v)$ has all components flagged
  ``trivial'', flag $v$ as ``trivial''.  Repeat the above procedure
  until no more flags are changed. Then flag all ``unknown yet'' words
  as ``trivial''.
\item[Computation] Let $w\in S^*$ be a word of which one asks whether
  it is trivial in $G$. If $w$ belongs to $V$, its flag answers the
  word problem. If $w$ acts non-trivially on $\Sigma$, it is
  non-trivial. Finally, if $w$ acts trivially on $\Sigma$, write
  $\psi(w)=(w_1,\dots,w_n)$. By property~(\ref{eq:contract}), each $w_i$
  is strictly shorter than $w$, so the algorithm can be applied
  inductively to it. $w$ is trivial if and only if all $w_i$ are
  trivial.
\end{description}
Only one point deserves a special justification, and that is the
flagging of ``unknown yet'' words as trivial. This is because such
words act trivially on the tree, so belong to
$\bigcap_{n\ge0}\stab_G(n)$, which by assumption is trivial.

\subsection{Proof of Theorem~\ref{thm:branch}}
\begin{lem}\label{lem:fg}
  Let $G=\langle S\rangle$ be a finitely generated group and let
  $H=\langle T\rangle$ be a finite-index subgroup of $G$, for some
  $T\subset S^*$. Let $\tilde S$ be a set in bijection with $S$, and
  for $w=s_1\dots s_n\in S^*$ set $\tilde w=\tilde s_1\dots\tilde
  s_n\in\tilde S^*$.
  
  There exists a finitely presented group $\Gamma=\langle \tilde
  S\rangle$ such that $\pi:\Gamma\overset{\tilde s\mapsto
    s}\longrightarrow G$ is an epimorphism, and $\pi^{-1}(H) = \langle
  \tilde T\rangle$ in $\Gamma$.
\end{lem}
\begin{proof}
  Consider first $F_S$ with its natural projection $\pi:F_S\to G$, and
  set $\Delta=\pi^{-1}(H)$. Since $\Delta$ has finite index in $F_S$,
  it is finitely generated, say by the set $U$. Our purpose is to find
  a quotient of $F_S$ in which $\Delta$ is generated by $T$. For each
  $u\in U$, let $w_u$ be an expression of $\pi(u)$ over $T$. It then
  suffices to consider
  \[\Gamma = \PRES{S}{u^{-1}w_u\,\forall u\in U}.\]
\end{proof}
In words, $\Gamma$ a finitely presented group such that the subgroup
lattice between $G$ and $\langle T\rangle$ is isomorphic to the
lattice between $\Gamma$ and $\langle T\rangle$, where the different
$\langle T\rangle$'s lie in different groups.

\begin{proof}[Proof of Theorem~\ref{thm:branch}]
  Let $G$ be regular branch on its subgroup $K$, and fix generating
  sets $S$ for $G$ and $T$ for $K$. It loses no generality to assume
  $K\le\stab_G(1)$, since one may always replace $K$ by
  $K\cap\stab_G(1)$. Let $\Gamma_0$ be the group given by
  Lemma~\ref{lem:fg} for $H=K$. Let $\Delta_0$ be the natural lift of
  $\stab_G(1)$ to $\Gamma_0$; and let $\Upsilon_0$ be the natural lift of
  $K$ to $\Gamma_0$.
  
  Let $U$ be a generating set of $\stab_1(G)$ (so $\Delta_0$ is
  generated by $\tilde U$), and let $\tilde\psi:\Delta_0\to\Gamma_0^d$
  be the natural lift of $\psi:\stab_G(1)\to G^d$ ; it maps $\tilde u$
  to $\widetilde{\psi(u)}$, where the wide tilde is applied to all $d$
  factors of $\psi(u)$. Note that $\tilde\psi$ satisfies the contracting
  condition for the same constant $D$ as $\psi$.
  
  Since $G$ is regular branch, there is an embedding $1^i\times
  K\times 1^{d-1-i}\hookrightarrow K$, from which for each generator
  $t\in T$ of $K$ we may choose a word $\phi_i(t)\in T^*$ such that
  $\psi(\phi_i(t))=(1,\dots,t,\dots,1)$ with the `$t$' in position
  $i$.
  
  Now write
  $\tilde\psi(\phi_i(t))=(r_{t,1},\dots,r_{t,i}t,\dots,r_{t,d})$ for
  some $r_{t,i}\in\Upsilon_0$. These elements' images in $K$ are
  trivial, since $\tilde\psi$ is a lift of $\psi$. Furthermore, since
  $\tilde\psi$ is contracting, one may replace $\{r_{t,i}\}_{t\in
  T,i\in\Sigma}$ by its iterates under all $\pi_i\tilde\psi$, where
  $\pi_i$ is the projection on the $i$-th factor, and still obtain a
  finite set of relations.

  Let $\Gamma$ be the quotient of $\Gamma_0$ by this sets' normal
  closure. Then $\Gamma$ is finitely presented and surjects onto $G$
  (since $r_{t,i}\cong1$ in $G$). Let $\Delta$ and $\Upsilon$ be the
  images of $\Delta_0$ and $\Upsilon_0$ in $\Gamma$, and note that
  $\psi$ lifts again to $\tilde\psi$ on $\Gamma$, because the new
  relators $r_{t_i}$ map to other new relators.

  The data are summed up in the following diagram, which should be
  viewed as a ``chair with $\psi$ and $\tilde\psi$ coming forward'':
  \[\begin{diagram}
    \node[2]{\Gamma}\arrow{s,-}\arrow{e,A}\node{G}\arrow{s,-}\\
    \node{\Gamma^d}\arrow{s,-}\node{\Delta}\arrow{s,-}
    \arrow{w,t}{\tilde\psi}\arrow{e,A}\node{\stab_G(1)}\arrow{s,-}
    \arrow{e,t,J}{\psi}\node{G^d}\arrow{s,-}\\
    \node{\Upsilon^d}\node{\Upsilon}\arrow{e,A}\node{K}\node{K^d}\arrow{w,L}
  \end{diagram}\]
  
  Since $\im\tilde\psi$ contains $\Upsilon^d$, it has finite index in
  $\Gamma^d$. Since $\Gamma^d$ is finitely presented, $\im\tilde\psi$
  too is finitely presented. Similarly, $\Delta$ is finitely
  presented, and we may express $\ker\tilde\psi$ as the normal closure
  $\langle R_1\rangle^\#$ in $\Delta$ of those relators in
  $\im\tilde\psi$ that are not relators in $\Delta$. Clearly $R_1$ may
  be chosen to be finite.

  We now use the assumption that $G$ is contracting, with constant $D$.
  Let $R_2$ be the set of words over $S$ of length at most $D$ that
  represent the identity in $G$. Set $R=R_1\cup R_2$, which clearly is
  finite.
  
  We now consider $T$ as a set distinct from $S$, and not as a subset
  of $S^*$. We extend each $\phi_i$ to a monoid homomorphism $(S\cup
  T)^*\to(S\cup T)^*$ by defining it arbitrarily on $S$.

  Assume $\Gamma=\SPRES SQ$, and let $w_t\in S^*$ be a
  representation of $t\in T$ as a word in $S$. We claim that the
  following is an \lpresn\ of $G$:
  \begin{equation}\label{eq:branchpres}
    G = \LPRES{S\cup T}{Q\cup\{t^{-1}w_t\}_{t\in T}}{%
      \{\phi_i\}_{i\in\Sigma}}{R_1\cup R_2}.
  \end{equation}
  For this purpose, consider the following subgroups $\Xi_n$ of
  $\Gamma$: first $\Xi_0=\{1\}$, and by induction
  \[\Xi_{n+1} = \big\{\gamma\in\Delta\big|\,
  \tilde\psi(\gamma)\in\Xi_n^d\big\}.\]
  We computed $\Xi_1=\langle R\rangle^\#$. Since $G$ acts transitively
  on the $n$-th level of the tree, a set of normal generators for
  $\Xi_n$ is given by $\bigcup_{i\le n}\phi^i(R)$, where
  $\phi$ is any choice of $\phi_i$ for $i\in\Sigma$. We also note that 
  $\psi(\Xi_{n+1})=\Xi_n^d$.

  We will have proven the claim if we show
  $G=\Gamma\big/\bigcup_{n\ge0}\Xi_n$. Let then $w\in\Gamma$ represent
  the identity in $G$. Applying to it $|w|$ times the map $\psi$, we
  obtain $d^{|w|}$ words that are all of length at most $K$, that is,
  that belong to $\Xi_1$. Then since $\psi(\Xi_{n+1})=\Xi_n^d$, we
  get $w\in\Xi_{|w|+1}$, and~(\ref{eq:branchpres}) is a presentation
  of $G$.
  
  As a bonus, the presentation~(\ref{eq:branchpres}) expresses $K$ as
  the subgroup of $G$ generated by $T$.
\end{proof}

A few remarks are in order. First, one can usually do with only one
substitution, say $\phi_1$, since in many cases the other $\phi_i$ are
conjugates of $\phi_1$. Second, $\phi_1$ induces an isomorphism from
$K$ to its subgroup $K\times1^{d-1}$, so there is a one-step \HNN\
extension of $G$ that is finitely presented --- namely, the extension
identifying $K$ and $K\times1^{d-1}$. Third, in many cases (but not
all) $\phi_1$ can be extended to an endomorphism of $G$; in that case,
one may delete $T$ from the generating set and obtain an ascending
\lpresn.

In all cases, $K$ admits an ascending \lpresn, so embeds in a finitely
presented group $L$, and $\langle G,L\rangle$ is a finite extension of
$L$, hence is a finitely presented group containing $G$.

\begin{proof}[Proof of Theorem~\ref{thm:notfp}]
  Since $G$ is contracting, there is a constant $D$ such that
  $|w_i|<|w|$ whenever $|w|>D$. This implies, using the triangular
  inequality, that there are constants $\eta<1$ and $D'$ such that
  $|w_i|<\eta|w|$ whenever $|w|>D'$.

  Now levels can be ``collapsed'' in a branch group: for any $k$ we
  may consider the (same) action of $G$ on $(\Sigma^k)^*$, with map
  $\psi$ given by $k$-fold composition of the original map $\psi$. The
  resulting group action is still branch.

  However, the result of this process is that the constant $\eta$
  above can be replaced by any power of itself, say $\frac12$, at the cost of
  enlarging the branching number of the tree.

  The generating set then now be replaced by a ball of sufficiently
  large radius, so that the constant $L$ becomes $1$.

  We have reached a ``canonical situation'', where the maps $\psi$
  and $\tilde\psi$ satisfy $|w_i|\le\frac12(|w|+1)$ for all $w$.
  
  Assume now by contradiction that $G$ is finitely presented, say
  $G=\SPRES SR$ with $\pi:F_S\to G$ the canonical map, and assume that
  the greatest length among the relators is minimal.  All
  $r\in R$ being trivial in $G$, satisfy \emph{a fortiori}
  $\pi(r)\in\stab_G(1)$, so $\tilde\psi(r)=(r_1,\dots,r_d)$ is well
  defined. By the Reidemeister-Schreier process, a presentation of
  $G\times1\times\dots\times1$ is
  \[\PRES S{r_i\text{ for all }r\in R\text{ and }i\in\Sigma}.\]
  By our assumptions that $|r_i|\le\frac12(|r|+1)$ and $\max|r|$
  is minimal, we must have $|r|\le1$ for all relations, so $G$ is
  free. However a free group may not contain commuting subgroups with
  trivial intersection, like $K\times1\times\dots\times1$ and
  $1\times\dots\times1\times K$. This is our required contradiction.
\end{proof}

\subsection{Schur multipliers}\label{subs:br:schur}
In his paper~\cite{grigorchuk:bath} Rostislav Grigorchuk computed
explicitly the Schur multiplier $H_2(G,\Z)$ of his group --- he proved
that $H_2(G)=(\Z/2)^\infty$. We outline here a general computation
$H_2(G,\Z)$ for branch groups $G$.

\begin{thm}
  Let $G$ be a finitely generated, contracting, semi-fractal, regular
  branch group. Then $H_2(G,\Z)\cong A\oplus B^\infty$, for finite
  abelian groups $A,B$.
\end{thm}
As a consequence, all such groups are infinitely presented.
\begin{proof}
  We concentrate on the exact sequence $1\to K^d\to K\to Q\to 1$, for
  some finite group $Q$. By~(\ref{eq:shapiro}) and~(\ref{eq:product}),
  $H_2(K,\Z Q)=H_2(K^d,\Z)=H_2(K,\Z)^d\oplus (K/K'\otimes
  K/K')^{d(d-1)/2}$. Taking $Q$-invariants of the right-hand side
  collapses all $d$ copies of $H_2(K,\Z)$ together, but we are left
  with the equation $H_2(K,\Z)=H_2(K,\Z)\oplus B$, where $B$, a
  quotient of $(K/K'\otimes K/K')^{d(d-1)/2}$, is a finite group.

  Then $H_2(G,\Z)$ is obtained from $H_2(K,\Z)$ by extension and
  quotient by finite-rank abelian groups, and the claimed result
  follows.
\end{proof}

Note, as a corollary, that if $K$ is perfect, then it is a finitely
presented, infinitely related group with trivial Schur multiplier.

\subsection{Perfect branch groups}\label{subs:perfect}
We consider a class of branch groups, of special interest for being
perfect. They form a subclass of the \textsf{GGS} groups studied
in~\cite{bartholdi-s:wpg}. Let $A$ be a finite, perfect, group acting
transitively on $\Sigma$, with two elements $\ast\neq\dagger\in\Sigma$
such that $\stab_A(\ast)\setminus\stab_A(\dagger)\neq\emptyset$ (think
for instance $\mathfrak A_5$).

Let $\overline A$ be a copy of $A$, and consider $\Gamma=A*\overline
A$. Define an action of $\Gamma$ on $\Sigma^*$ by
\begin{align*}
  (\sigma_1\sigma_2\dots\sigma_n)^a&=\sigma_1^a\sigma_2\dots\sigma_n,\\
  (\sigma_1\sigma_2\dots\sigma_n)^{\overline a}&=\begin{cases}
    \sigma_1\sigma_2^a\sigma_3\dots\sigma_n&\text{ if }\sigma_1=\ast,\\
    \sigma_1(\sigma_2\sigma_3\dots\sigma_n)^{\overline a}&\text{ if }\sigma_1=\dagger,\\
    \sigma_1\sigma_2\dots\sigma_n&\text{ otherwise}.
  \end{cases}
\end{align*}
Let $G$ be the group defined by this action, i.e.\ the quotient of
$\Gamma$ by the kernel of the action.

\begin{prop}
  $G$ is a perfect finitely generated regular branch group over
  itself.
\end{prop}
\begin{proof}
  Clearly $G$ is perfect, being generated by two perfect groups, and
  finitely generated, being generated by two finite groups.
  
  Note now that $\stab_G(1)=\overline A^G$. The map $\psi:G\to
  G\wr_\Sigma A$ is given by
  \[\psi(a)=(1,\dots,1)a,\qquad\psi(\overline a)=(a,1,\dots,1,\overline a)1,\]
  where in this last expression the `$a$' is at position $\ast$ and
  the `$\overline a$' is at position $\dagger$.  The conditions on $A$
  imply that it contains an element $x$ moving $\dagger$ but not
  $\ast$.  The computation $\psi[\overline a,\overline
  b^x]=([a,b],1,\dots,1)$ shows that $\psi(G)$ contains
  $A\times1\dots\times\dots\times1$, since $A$ is perfect; then
  $\psi(G)$ contains too
  \[(a,1,\dots,1)^{-1}\psi(\overline a)
  =(a^{-1},1,\dots,1)(a,1,\dots,1,\overline a)=(1,\dots,1,\overline
  a),\]
  so $\psi(G)$ contains $1\times\dots\times1\times\overline A$, and
  since $A$ is $Y$-transitive it contains $G\times\dots\times G$.
  (Explicitly, we have $\psi^{-1}(G^\Sigma)=\stab_G(1)$.)
\end{proof}

In this context, the statements of the previous section simplify: we
have perfect regular branch groups $G$ with $H_2(G,\Z)=0$, that are
finitely \lpresd\ but infinitely presented. The group
$\Gamma=A*\overline A$ is the same as the $\Gamma$ of
Theorem~\ref{thm:branch}, and the subgroup $\Delta$ is {\LARGE $*$}${}_{a\in
  A}\overline A^a$. The map $\tilde\psi$ is given by
\[\tilde\psi(\overline a^b) = (a\text{ in position }\ast^b,\overline
a\text{ in position }\dagger^b).\]
Let $\phi$ be some word substitution on $\Gamma$ mapping $g$ to
$(g,1,\dots,1)$, as given by the computations in the previous theorem.
We then have an \lpresn
\[G = \left\langle A,\overline A\left\|\phi\left|
      \begin{array}{ll}
        \overline a^{1-b}&\text{ whenever }\ast^b=\ast,\dagger^b=\dagger\\
        \overline a^{1-b+c-d}&\text{ whenever }\ast^b=\ast,\dagger^b=
        \dagger^c,\ast^c=\ast^d,\dagger^d=\dagger\text{ are all distinct}\\
        \overline a^{1-b+c-d(1-b+c)}&\text{ whenever }\ast^b=\dagger^c=
        \dagger^d=\ast,\dagger^b=\ast^c,\ast^d=\dagger\text{ are all distinct}\\{}
        [\phi(a),\phi(b)^c]&\text{ whenever }\ast^c\neq\ast
      \end{array}\right.\right.\right\rangle.\]
Indeed, the first three relations identify all products $\overline
a^{b_1+\dots+b_n}$ with same $\psi$-image, and the last ones are the
commutation relations lifted from $\Gamma\times\dots\times\Gamma$.

\section{Examples}\label{sec:examples}
We start by an example of wreath product:
\begin{thm}\label{thm:ll}
  The following is an \lpresn\ of the ``lamplighter group''
  $G=\Z/2\wr\Z$:
  \[G = \LPRES{a,b,t}{a^2,a^{-1}b}{\phi}{[a,b]},\]
  where $\phi:\{a,b,t\}^*\to\{a,b,t\}^*$ is given by
  \[\phi(a)=a,\qquad\phi(b)=b^t,\qquad\phi(t)=t.\]

  However, this group admits no finite presentation.
\end{thm}
\begin{proof}
  A presentation of $G$ is
  \[\PRES{a,t}{a^2,[a,a^{t^i}]\,\forall i\in\Z}.\]
  By conjugating the last relation by $t^i$, we may assume
  the set of relators is $a^2$ and the $[a,a^{t^i}]$ with $i\ge0$. The
  latter are precisely the relators obtained from $\phi^i([a,b])$ by
  applying the relation $a=b$.
  
  It follows from~\cite{baumslag:fpwreath} that $G$ is not finitely
  presented. Even better, (\ref{eq:blackburn}) gives
  $H_2(\Z/2\wr\Z,\Z)=(\Z/2)^\infty$.
\end{proof}

Note however the following seemingly similar example, due to Gilbert
Baumslag, which is finitely \lpresd\ by arguments similar to those in
Theorem~\ref{thm:ll}:
\begin{thm}[\cite{baumslag:trivialsm}]\label{thm:baumslag}
  The group
  \[G = \PRES{a,b,t}{a^{t-4},b^{2t-1},[a,b^{t^i}]\,\forall i\in\Z}\]
  is an infinitely-presented, metabelian group, with $H_2(G,\Z)=0$.
\end{thm}
This example was devised to show that the Schur multiplier's rank may
be much smaller than the number of relators. In that view, we may ask
the following question:
\begin{problem}
  Do there exist non-finitely-\lpresd\ groups with trivial Schur
  multiplier?
\end{problem}

An interesting example of group acting on a rooted tree is the
``Brunner-Sidki-Vieira group''; we rephrase their result in terms of
$L$-presentations:
\begin{thm}[\cite{brunner-s-v:nonsolvable}, Proposition 15]
  Consider the group $G=\langle\mu,\tau\rangle$ acting on $\{1,2\}^*$,
  with $\psi(a^{-1}\mu)=(1,\mu^{-1})$ and $\psi(a^{-1}\tau)=(1,\tau)$
  (so $\tau$ and $\mu$ act like $a$ on the top node of the tree. Note
  that $G$ is neither rooted nor branch, though it is ``weakly
  branch''~\cite{grigorchuk:jibg}.) Writing $\lambda=\tau\mu^{-1}$,
  $G$ admits the ascending $L$-presentation
  \[G = \ELPRES{\lambda,\tau}{\phi}{[\lambda,\lambda^\tau],
  [\lambda,\lambda^{\tau^3}]},\]
  where $\phi$ is defined by $\tau\mapsto\tau^2$ and
  $\lambda\mapsto\tau^2\lambda^{-1}\tau^2$.
\end{thm}
We may even conclude that $H_2(G,\Z)=(\Z\times\Z)^\infty$, freely
generated by the images of $\phi^n[\lambda,\lambda^\tau]$ and
$\phi^n[\lambda,\lambda^{\tau^3}]$.

We now give presentations for four of the five ``testbed groups''
studied in~\cite{bartholdi-g:cras-spectrum,bartholdi-g:parabolic}.
\subsection{An \lpresn\ for $G$}
The group $G$, the \emph{first Grigorchuk group}, is the $2$-rooted
group $G=\langle a,b,c,d\rangle$, with $a$ the rooted element and
$b,c,d$ defined by
\[\psi(b)=(a,c),\quad\psi(c)=(a,d),\quad\psi(d)=(1,b).\]
$G$ is a regular branch group over $K=\langle(ab)^2\rangle^\#$.

\begin{thm}\label{thm:G}
  The Grigorchuk group $G$ admits the ascending \lpresn
  \[G = \ELPRES{a,c,d}{\sigma}{a^2,[d,d^a],[d^{ac},(d^{ac})^a]},\]
  where $\sigma:\{a,c,d\}^*\to\{a,c,d\}^*$ is defined by
  \[\sigma(a)=aca,\quad\sigma(c)=cd,\quad\sigma(d)=c.\]
\end{thm}
\begin{proof}
  Rephrasing of~\cite{lysionok:pres}.
\end{proof}

\subsection{An \lpresn\ for $\tilde G$}
The group $\tilde G$, the \emph{Grigorchuk supergroup}, is the
$2$-rooted group $G=\langle a,\tilde b,\tilde c,\tilde d\rangle$, with
$a$ the rooted element and $\tilde b,\tilde c,\tilde d$ defined by
\[\psi(\tilde b)=(a,\tilde c),\quad\psi(\tilde c)=(1,\tilde d),
\quad\psi(\tilde d)=(1,\tilde b).\]
$\tilde G$ is a regular branch group over $\tilde K=\langle(a\tilde
b)^2,(a\tilde d)^2\rangle^\#$. It is named thus because it contains
$G$ as a subgroup.

\begin{thm}\label{thm:Gt}
  The group $\tilde G$ admits the ascending \lpresn
  \[\tilde G=\ELPRES{a,\tilde b,\tilde c,\tilde d}{%
  \tilde\sigma}{a^2,[\tilde b,\tilde c],
  [\tilde c,\tilde c^a],[\tilde c,\tilde d^a],[\tilde d,\tilde d^a],
  [\tilde c^{a\tilde b},(\tilde c^{a\tilde b})^a],
  [\tilde c^{a\tilde b},(\tilde d^{a\tilde b})^a],
  [\tilde d^{a\tilde b},(\tilde d^{a\tilde b})^a]},\]

  where $\tilde\sigma:\{a,\tilde b,\tilde c,\tilde d\}^*\to\{a,\tilde
  b,\tilde c,\tilde d\}^*$ is defined by
  \[a\mapsto a\tilde ba,\quad\tilde b\mapsto\tilde d,\quad\tilde
  c\mapsto\tilde b,\quad\tilde d\mapsto\tilde c.\]
\end{thm}
\begin{proof}
  Rephrasing of~\cite[Proposition~5.6]{bartholdi-g:parabolic}.
\end{proof}

\subsection{An \lpresn\ for $\Gamma$}
The group $\Gamma$, the \emph{Fabrykowski-Gupta group}, is the
$3$-rooted group $G=\langle a,r\rangle$, with $a$ the rooted element
and $r$ defined by
\[\psi(r)=(a,1,r).\]
$\Gamma$ is a regular branch group over $\Gamma'=\langle[a,r]\rangle^\#$.

\begin{thm}\label{thm:Gamma}
  The Fabrykowski-Gupta group $\Gamma$ admits the ascending
  \lpresn
  \[\ELPRES{a,r}{\sigma,\chi_1,\chi_2}{a^3,[r^{1+a^{-1}-1+a+1},a],
    [a^{-1},r^{1+a+a^{-1}}][r^{a+1+a^{-1}},a]},\]
  where $\sigma,\chi_1,\chi_2:\{a,r\}^*\to\{a,r\}^*$ are given by
  \begin{xalignat*}{2}
    \sigma(a)&=r^{a^{-1}},&\sigma(r)&=r,\\
    \chi_1(a)&=a,&\chi_1(r)&=r^{-1},\\
    \chi_2(a)&=a^{-1},&\chi_2(r)&=r.
  \end{xalignat*}
\end{thm}

\begin{proof}
  We follow Theorem~\ref{thm:branch}. Consider first the group
  $F=\SPRES{a,r}{a^3,r^3}$. Clearly,
  $F/F'\cong(\Z/3)^2\cong\Gamma/\Gamma'$. Using the computer algebra
  program \textsc{Gap}, we compute a presentation for $\im\tilde\psi$,
  and rewrite its relators as words in $X$, where $X$ is a generating
  set for $\Gamma'$. We also construct a group homomorphism
  $\sigma_0:\Gamma'\to1\times1\times\Gamma'$. Then
  Theorem~\ref{thm:branch} gives a finite \lpresn\ for
  $\Gamma$ with generators $\{a,r\}\cup X$.
  
  We now note that $\sigma_0$ can be extended to a homomorphism
  $\sigma:\Gamma\to A\times
  R\times\Gamma$, where $A=\langle a\rangle$ and $R=\langle
  r\rangle$ have order $3$.
  The substitution $\sigma$ can be used instead of $\sigma_0$, giving rise
  to a simpler presentation with generators $a,r$.

  Finally, we note that the presentation can be simplified from $6$
  iterated relators to $2$ by introducing two extra substitutions
  $\chi_1,\chi_2$ induced by group automorphisms.
\end{proof}

Note that for $G$ and $\tilde G$ the iterated relations are of the
form $[x,x^a]$ where $x$ belongs to a first-level rigid
stabilizer. For $\Gamma$, however, one obtains fewer relations by
considering more general expressions, as above.

\subsection{An \lpresn\ for $\overline\Gamma$}
The group $\overline\Gamma$, introduced
in~\cite{bartholdi-g:cras-spectrum}, is the $3$-rooted group $G=\langle
a,s\rangle$, with $a$ the rooted element and $s$ defined by
\[\psi(s)=(a,a,s).\]
Set $x=ta^{-1},y=a^{-1}t$ and $K=\langle x,y\rangle$, a torsion-free
index-$3$ subgroup of $\overline\Gamma$.

\begin{thm}\label{thm:GammaB}
  The groups $K$ and $\overline\Gamma$ are not branch, but are finitely
  \lpresd.
\end{thm}
\begin{proof}
  We start by computing an \lpresn\ for $K$. As above,
  $\psi(K')$ contains $K'\times K'\times K'$; but $K/K'\cong\Z^2$ and
  neither $\overline\Gamma$ nor $K$ are branch.

  First we chose generators of $\stab_K(1)$:
  \begin{xalignat*}{2}
    \alpha&=x^{-1}y=(x,1,x^{-1}),&
    \beta&=y^{-1}x^{-1}y^{-1}=(y,1,y^{-1}),\\
    \gamma&=x^{-1}y^{-1}x^{-1}=(1,x,x^{-1}),&
    \delta&=xy^{-1}=(1,y,y^{-1}).
  \end{xalignat*}

  Then, we chose generators of $K'$:
  \begin{xalignat*}{2}
    e&=\beta^{-1}\delta\gamma=[y,xyx]=(y^{-1},yx,x^{-1}),&
    f&=\gamma\beta^{-1}\delta=[x^{-1}y^{-1}x^{-1},y]=(y^{-1},xy,x^{-1}),\\
    g&=\gamma^{-1}\alpha\beta=[x,y]=(xy,x^{-1},y^{-1}),&
    h&=\beta\gamma^{-1}\alpha=[y^{-1}x^{-1},y^{-1}]=(yx,x^{-1},y^{-1}).
  \end{xalignat*}

  The fact that $K'$ is normal can be seen in the following
  conjugation relations:
  \[\begin{array}{c|cccc}
    \swarrow & e & f & g & h\\ \hline
    x & g^{-1}hf^{-1}g^{-1} & f^{-1}g^{-1} & e & g^{-1}he \\
    x^{-1} & ehf^{-1} & h^{-1}e^{-1} & h^{-1}f^{-1} \\
    y & g^{-1}e^{-1} & h^{-1}e^{-1} & g^{-1}eh & f \\
    y^{-1} & e^{-1}gf & h & f^{-1}g^{-1} & f^{-1}g^{-1}ef
  \end{array}\]

  Define now the group $L$ by its generators $S=\{x^{\pm1},y^{\pm1}\}$
  and relators $k^s=w_{k,s}$ for $s\in S$ and $k\in\{e,f,g,h\}$, where
  $w_{k,s}$ is the word in the above table. Note then that $L'$ is
  generated by the words $e,f,g,h$.
  
  As in the proof of Theorem~\ref{thm:branch}, consider the map
  $\tilde\psi:\stab_L(1))\to L^3$ corresponding to $\psi:\stab_K(1)\to
  K^3$, and by the Reidemeister-Schreier method compute a presentation
  for $\tilde\psi(\stab_L(1))$.
  
  Assume $L=\SPRES SQ$. A presentation for $L^3$ is
  $\SPRES{S_1\cup S_2\cup S_3}{Q_1\cup Q_2\cup Q_3\cup[S_i,S_{j\neq
    i}]}$. The image of $\tilde\psi$ can be described as
  $\{(u,v,w)\in L^3|\,uvw\in L'\}$. We choose
  $\{x_3^my_3^n\}_{m,n\in\Z}$ as Schreier transversal for this
  subgroup, and denote the Schreier generators
  \begin{xalignat*}{2}
    u_{imn} &= x_3^my_3^nx_iy_3^{-n}x_3^{-m-1}, &
    v_{imn} &= x_3^my_3^ny_iy_3^{-n-1}x_3^{-m}.
    \intertext{Then these generators satisfy}
    u_{1mn}&=u_{100}=\alpha,&
    v_{1mn}&=u_{200}^{-m}v_{100}u_{200}^m=\gamma^{-m}\beta\gamma^m,\\
    u_{2mn}&=u_{200}=\gamma,&
    v_{2mn}&=u_{100}^{-m}v_{200}u_{100}^m=\alpha^{-m}\delta\alpha^m,\\
    u_{3mn}&=u_{200}^{-m}v_{100}^{-n}u_{200}^{-1}v_{100}^nu_{200}^{m+1},&
    v_{3mn}&=1.
  \end{xalignat*}

  The relators we obtain, in terms of $\alpha,\beta,\gamma,\delta$,
  are
  \begin{align*}
    [\alpha,\gamma],\;[\alpha\beta,\gamma\delta],&\\
    [\alpha\gamma^{-1},\beta^{-n}\gamma^{-1}\beta^n],\;
    [\gamma^{-n}\beta\gamma^n,\alpha^{-n}\delta\alpha^n]&\text{ for all
      }n\in\Z,\\
    w(\alpha\gamma^{-1},\beta),\;
    w(\alpha^{-1}\gamma,\delta),\;
    w(\alpha^{-1},\delta^{-1})&\text { for all }w\in Q.
  \end{align*}
  
  This clearly gives a finite \lpresn\ for
  $\tilde\psi(\stab_L(1))$ --- compare with the proof of
  Theorem~\ref{thm:ll}. Now the computation of a presentation for $K$
  can be finished as in the proof of Theorem~\ref{thm:branch}.

  Finally, a finite \lpresn\ for $G$ can be obtained using
  Proposition~\ref{prop:ext}.
\end{proof}

\subsection{An \lpresn\ for $\doverline\Gamma$}
The group $\doverline\Gamma$, the \emph{Gupta-Sidki group}, is the
$3$-rooted group $G=\langle a,t\rangle$, with $a$ the rooted element
and $t$ defined by
\[\psi(t)=(a,a^{-1},t).\]
$\doverline\Gamma$ is a regular branch group over
$\doverline\Gamma'=\langle[a,t]\rangle^\#$.

\begin{thm}\label{thm:GammaBB}
  The Gupta-Sidki group $\doverline\Gamma$ admits the \lpresn
  \[\LPRES{a,t,u,v}{a^3,t^3,u^{-1}t^a,v^{-1}t^{a^{-1}}}{\sigma,\chi}{%
    (tuv)^3,[v,t][vt,u^{-1}tv^{-1}u],[t,u]^3[u,v]^3[t,v]^3},\]
  where $\sigma,\chi:\{t,u,v\}^*\to\{t,u,v\}^*$ are given by
  \begin{gather*}
    \sigma:\begin{cases}
      t\mapsto t,\\
      u\mapsto [u^{-1}t^{-1},t^{-1}v^{-1}]t=u^{-1}tv^{-1}tuvt^{-1},\\
      v\mapsto t[tv,ut] = t^{-1}vutv^{-1}tu^{-1},
    \end{cases}\qquad
    \chi:\begin{cases}
      t\mapsto t^{-1},\\ u\mapsto u^{-1},\\ v\mapsto v^{-1}.
    \end{cases}
  \end{gather*}
\end{thm}
Note that $\chi$ is induced by the automorphism of $\doverline\Gamma$
defined by $a\mapsto a,\;t\mapsto t^{-1}$; however, $\sigma$ does not
extend to an endomorphism of $\doverline\Gamma$.

Note also that all the iterated relators can be expressed as words
over $\{t,u,v\}$ with $0$-sum in each variable. Their most natural
representation is as closed paths in the $\{t,u,v\}$-space:
\[\begin{picture}(364,102)
  \put(0,8){\epsfig{file=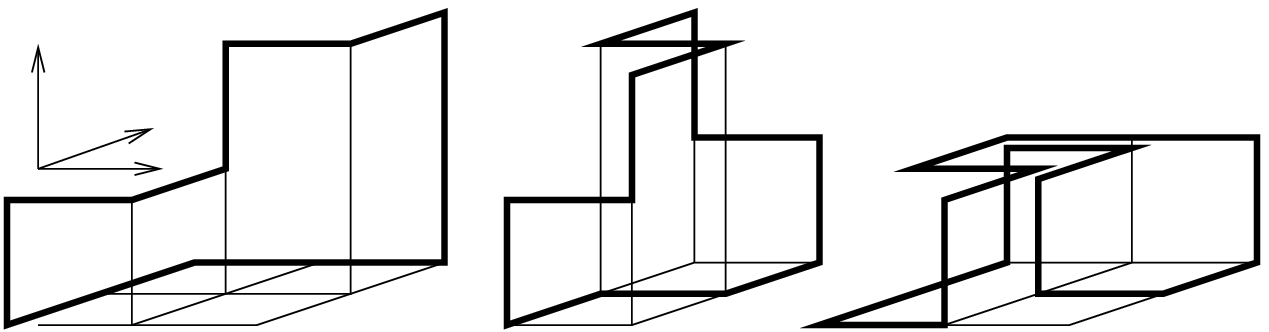}}
  \put(6,70){$t$}
  \put(26,50){$u$}
  \put(22,64){$v$}
  \put(50,0){$(tuv)^3$}
  \put(150,0){$[v,t][vt,u^{-1}tv^{-1}u]$}
  \put(270,0){$[t,u]^3[u,v]^3[t,v]^3$}
\end{picture}\]
Then the fact that these elements are non-trivial relations translates
to the fact that their projection on any plane $t=-u$, $u=-v$ or
$v=-t$ gives a trivial path (up to $t^3=u^3=v^3=1$), while they
themselves are not trivial paths. Incidentally, these projections are
none but the $\psi_i:\langle t,u,v\rangle\to\langle a,t\rangle$, for
$i\in\Sigma$.

\begin{proof}
  We follow Theorem~\ref{thm:branch}. Consider first the group
  $F=\SPRES{a,t}{a^3,t^3}$. Clearly,
  $F/F'\cong(\Z/3)^2\cong\doverline\Gamma/\doverline\Gamma'$. Using the
  computer algebra program \textsc{Gap}, we compute a presentation for
  $\im\tilde\psi$, and rewrite its relators as words in $X$, where $X$ is
  a generating set for $\doverline\Gamma'$. We also construct a group
  homomorphism
  $\sigma_0:\doverline\Gamma'\to1\times1\times\doverline\Gamma'$. Then
  Theorem~\ref{thm:branch} gives a finite \lpresn\ for
  $\doverline\Gamma$ with generators $\{a,t\}\cup X$.
  
  We now note that $\sigma_0$ can be extended to a homomorphism
  $\sigma:\stab_{\doverline\Gamma}(1)\to A\times
  A\times\doverline\Gamma$, where $A=\langle a\rangle$ has order $3$.
  The substitution $\sigma$ can be used instead of $\sigma_0$, giving rise
  to a simpler presentation with generators $a,t,u,v$, where
  $t,u=t^a,v=t^{a^{-1}}$ is a generating set for
  $\stab_{\doverline\Gamma}(1)$.

  Finally, we note that the presentation can be simplified from $6$
  iterated relators to $3$ by introducing an extra substitution $\chi$,
  induced by a group automorphism.
\end{proof}

\begin{problem}
  Does there exist a finite ascending \lpresn\ for $\doverline\Gamma$?
\end{problem}

In these examples, easy computations yield $H_2(G,\Z)=H_2(\tilde
G,\Z)=(\Z/2)^\infty$ and
$H_2(\Gamma,\Z)=H_2(\doverline\Gamma,\Z)=(\Z/3)^\infty$.

%
%
%
%
%

\def\isbn#1{#1}
\bibliography{mrabbrev,people,math,grigorchuk,bartholdi}
\end{document}